\documentclass[11pt]{amsart}
\usepackage{geometry,amsrefs}                
\usepackage{graphicx}
\usepackage[usenames,dvipsnames]{color}
\usepackage{amssymb}
\usepackage{epstopdf}
\newtheorem{prop}{Proposition}
\newtheorem{lem}{Lemma}
\newtheorem{theorem}{Theorem}
\newtheorem{remark}{Remark}
\newtheorem{cor}{Corollary}
\title{Singular multicontact structures}
\author{Alessandro Ottazzi, Gerd Schmalz}

\begin{document}

\thanks{The authors have been partially supported by the Gruppo Nazionale per l'Analisi Matematica, la Probabilit\`a e le loro Applicazioni (GNAMPA) of the Istituto Nazionale di Alta Matematica (INdAM) and the ARC Discovery grant DP130103485}
\maketitle

\begin{abstract}
We describe the automorphisms of a singular multicontact structure, that is a generalisation of the Martinet distribution. 
Such a structure is interpreted  as a para-CR structure on a hypersurface $M$ of a direct product space $\mathbb R_+^2\times \mathbb R_-^2$. We introduce the notion of a finite type singularity analogous to CR geometry and, along the way, we prove extension results for para-CR functions and mappings on embedded para-CR manifolds into the ambient space. 
\end{abstract}

\section{Introduction}
Multicontact structures have been studied under different names and in diverse contexts since the end of the 19th century. 
Recently they have been treated in great generality in the context of  parabolic geometry and sub-Riemannian geometry \cite{CDeKR1, CDeKR2, KorMC, Ottazzi-Hess, Han-Oh-Schmalz, CCDEM, DeOtt}, see also \cite{Yatsui1, Yatsui2}.

The simplest instance of a multicontact structure is a 3-dimensional manifold with two direction fields $X,Y$ that span a contact distribution. This structure appears in the point-transformation geometry of 2nd order ODE and was first studied by Lie \cites{LieIII, LieIV, Lie1891} and Tresse \cite{Tresse}.
It is well-known that this multicontact geometry is analogous to the Cauchy-Riemann geometry of real hypersurfaces in $\mathbb C^2$, which was one of Cartan's \cites{Car33, Car32} inspirations for developing his technique of moving frames. The analogy between the two geometries has been described by Nuwrowski and Sparling \cite{MR2024797} in the intrinsic setup, but it seems that the extrinsic approach due to Chern and Moser \cite{CM} has never been applied to multicontact structures. For this approach the ambient complex manifold needs to be replaced by a direct product manifold. The analogy between complex manifolds and direct product manifolds has been summarised in the recent survey article \cite{HL} by Harvey and Lawson.

According to this approach one would define a para-CR structure on a hypersurface $M$ of a direct product space $\mathbb R_+^2\times \mathbb R_-^2$ (or more in general $\mathbb R_+^n\times \mathbb R_-^n$) as the structure on $TM$ induced by the embedding, namely the two direction fields (or more in general $n-1$-dimensional distributions) $TM\cap T\mathbb R_\pm^2$. In order to make this a multicontact structure we need to impose the condition that the commutator of the two direction fields generates the missing direction in $TM$ at each point. This is analogous to Levi non-degeneracy of a CR manifold.

In this paper we use an extrinsic approach to study singular multicontact structures, i.e. multicontact structures for which the contact condition fails on a thin subset.
In the context of sub-Riemannian geometry, the structures that we study define singular sub-Riemannian structures, and are a generalization of 
 the Martinet distribution \cite{Martinet, Cal-Cha-East, Zel-Zhit}. 
 
 We show that, in contrast to the CR situation, any para-CR structure can be locally realised by an embedding.
We define a notion of finite type singularities that is analogous to the notion of finite type in CR-geometry, and we study the homogeneous models and their symmetries. 

The article is structured as follows: In Section 2 we prove that any abstract para-CR manifold $M$ can be locally realised as an embedded para-CR manifold. Section 3 is devoted to the extension of para-CR functions to para-holomorphic functions in a neighbourhood of an embedded para-CR manifold. In Section 4 we define the notion of para-CR mappings, para-CR automorphisms and infinitesimal para-CR automorphisms and we show that they extend as para-holomorphic objects in a neighbourhood of an embedded para-CR manifold. In Section 5 we introduce the notion of finite type para-CR structures, which is analogous to CR-manifolds of finite type and derive a special form of the defining equation of an embedded para-CR manifold of finite type $k$. Section 6 is the central part of this article. There we compute the symmetries of the homogeneous models. 

\section{Embedding of para-CR structures}\label{emsec}
Let $M$ be a hypersurface in $\mathbb R_{xy}^2\times \mathbb R_{ab}^2$, locally given as 
$$y=a+ \phi(a,b,x)$$
where $\frac{\partial \phi}{\partial a}(0)=0$. This embedding  distinguishes two direction fields $TM\cap T\mathbb R_{xy}^2$ and $TM\cap T\mathbb R_{ab}^2$ on $M$, which in local coordinates $x,a,b$ take the form
$$X=\partial_x, \quad Y= (1+ \frac{\partial \phi}{\partial a})\partial_b - \frac{\partial \phi}{\partial b}\partial_a.$$

 We will call a 3-dimensional manifold $M$ with two distinguished direction fields $X,Y$ a para-CR manifold. If the distribution spanned by $X,Y$ is contact, i.e., if $[X,Y]$ generates the missing direction then we call the structure a multicontact structure. We use the term singular multicontact structure if the contact condition is not satisfied on some submanifold of $M$.

We show now that any abstract para-CR manifold can be locally embedded in this way. 

\begin{lem}\label{para-CR-embed}
Let $X$, $Y$ be two linearly independent vector fields in some neighborhood of the origin in $\mathbb R^3$. Then there exist coordinates $(x,a,b)$ such that $X=\mu\partial_x$ and $Y=\lambda(\partial_b +\psi(x,a,b)\partial_a)$,  where $\mu, \lambda$ are some non-vanishing functions.
\end{lem}

Proof. Let $(x^*,a^*,b^*)$ be coordinates such that $X=\partial_{x^*}$ and $Y(0)=\partial_{b^*}|_0$. Then $Y=\lambda(\partial_{b^*}+A\partial_{a^*}+B\partial_{x^*})$, where $\lambda, A,B$ are functions with $\lambda(0)=1$.

We apply a coordinate change 
$$a=a^*, \qquad b=b^*, \qquad x= \chi(x^*,a^*,b^*).
$$ 
Now, $X=\frac{\partial \chi}{\partial x^*} \partial_x$ and the $\partial_x$ component of $Y$ is
$$\lambda(\frac{\partial \chi}{\partial b^*}+ A \frac{\partial \chi}{\partial a^*} +B\frac{\partial \chi}{\partial x^*}).$$
It remains to choose $\chi$ as a solution of the PDE
$$\frac{\partial \chi}{\partial b^*}+ A \frac{\partial \chi}{\partial a^*} +B\frac{\partial \chi}{\partial x^*}=0. \quad\Box$$\medskip

\begin{prop}
Let $M$ be a $3$-manifold with a para-CR structure $(X,Y)$. Then $M$ can be locally embedded into $\mathbb R^4$ in such a way  that the induced para-CR structure coincides with the original one. 
\end{prop}

Proof. By Lemma \ref{para-CR-embed}, we may assume that $X=\partial_x$ and $Y=\partial_b+\psi \partial_a$ defines the para-CR structure on a coordinate chart in $\mathbb R^3$ with coordinates $x,a,b$.

We embed $M$  into $\mathbb R^4$ with coordinates $(a,b,x,y)$ as a hypersurface 
\begin{equation}\label{embed}
S\colon y=a+ \phi(a,b,x)=\tilde{\phi}(a,x,b)
\end{equation}
with $\frac{\partial \phi}{\partial a}|_{b=0}=0$.

Clearly, $X=\partial_x$ lifts to $S$ as $\partial_x + \frac{\partial \phi}{\partial x}\partial_y$, which is a section of $T\mathbb R^2_{xy}$. If we choose $\tilde{\phi}$ to be the solution of the Cauchy problem
\begin{align*}
\frac{\partial \tilde{\phi}}{\partial b} + \psi(a,b,x) \frac{\partial \tilde{\phi}}{\partial a}&=0\\
\tilde{\phi}|_{b=0}&=a
\end{align*}
then 
$$\partial_b - \frac{\frac{\partial \phi}{\partial b}}{1+\frac{\partial \phi}{\partial a} }\partial_a$$
is the direction field $TS\cap T\mathbb R^2_{ab}$. \hfill $\Box$

%

\section{para-CR functions}
We call a function $f=(u,v)\colon D\to \mathbb R^2$ from a domain $D\subset \mathbb R^{2n}$ with coordinates $(x_1,\dots,x_n,a_1,\dots,a_n)$ para-holomorphic if $\frac{\partial u}{\partial a_j}=0$ and $\frac{\partial v}{\partial x_j}=0$ for all $j=1,\dots n$.

Let $(M,X,Y)$ be a para-CR manifold. A function $f=(u,v)\colon M \to \mathbb R^2$ is called a para-CR function if $Xv=Yu=0$.


\begin{prop}\label{paraCRfunc}
If $(x,a,b)$ are local canonical coordinates on $M$ with $X=\partial_x$ and $Y=\partial_b+\psi\partial_a$, then $f=(u,v)$ is para-CR if and only if $$u=u(x,y), \qquad v=v(a,b)$$
where $y=a+\phi(x,a,b)$ is the function from the embedding \eqref{embed}. 
\end{prop}

Proof. It is easy to see that $Xv=0$ is equivalent to $v=v(a,b)$. After a change of coordinates to
$x^*=x$, $b^*=b$, $y=a+\phi(x,a,b)$ the vector field $Y$ becomes $\partial_b$ and $Yu=0$ is now equivalent to $u=u(x^*,y)=u(x,y)$. \hfill $\Box$

\begin{prop}
Suppose that $M$ is embedded into $\mathbb R^4$ as described in Section \ref{emsec}. Then a para-CR function $f=(u,v)$ extends in a unique way to a para-holomorphic function $\tilde{f}=(\tilde{u},\tilde{v})$ in some neighbourhood of $M$. 
\end{prop}

Proof. Let $(x,y,a,b)$ be the coordinates of the embedding, as above. According to the characterisation of para-CR functions in Proposition \ref{paraCRfunc}, $u=u(a,b)$ extends to all $x,y$ as $\tilde{u}(x,y,a,b)=u(a,b)$ and  $v=v(x,y)$ extends to all $a,b$ as $\tilde{v}(x,y,a,b)=v(x,y)$. \hfill $\Box$

Notice  that, in contrast to the CR case, we have a two-sided extension and no assumption on the Levi form is made. 

\section{Extension of para-CR mappings}
Let $(M,X,Y)$ and $(M',X',Y')$ be two para-CR manifolds. A differentiable mapping $F\colon M\to M'$ is called a para-CR mapping if $F^*X'=\lambda X$ and $F^*Y'=\mu Y$ for some non-vanishing functions $\lambda, \mu$.

\begin{prop}\label{extmap}
If $\phi\colon M\to M'$ is a para-CR mapping of embedded para-CR manifolds with embeddings $(x,y,a,b)$ and $(x',y',a',b')$ then $F$ extends to a para-holomorphic mapping $\tilde{F}$ from a neighbourhood of $M$ to a neighbourhood of $M'$. 
\end{prop}

Proof. If $M'$ is embedded into $\mathbb R^4$ by $(x',y',a',b')$ then the composition of $F$ with the coordinate maps gives two para-CR functions $f_1=(x',a')\circ F$ and $f_2=(y',b')\circ F$ on $M$. If $M$ is embedded into  $\mathbb R^4$ by $(x,y,a,b)$ then $f_1$ and $f_2$ extend as para-holomorphic functions $\tilde{f}_1,\tilde{f}_2$ to a neighbourhood of $M$. Now $\tilde{F}=(\tilde{f}_1,\tilde{f}_2)$.\hfill $\Box$\bigskip

We define a local para-CR automorphism of a para-CR manifold $M$ as a diffeomorphism $F$ that is at the same time a para-CR mapping from an open subset $U\subset M$ onto another open subset $U'\subset M$. A locally given vector field $V$ is called a local infinitesimal para-CR automorphism of $M$ if the flow of $V$ consists of local para-CR automorphisms. Proposition \ref{extmap} readily implies the following corollary.
\begin{cor} Let $V$ be a local infinitesimal para-CR automorphism of an embedded para-CR manifold $M\subset \mathbb R^4$. Then $V$ extends to a para-holomorphic vector field $\tilde{V}$ in some neighbourhood of $M$. Conversely, the restriction of any para-holomorphic vector field $\tilde{V}$, such that $\tilde{V}|_M$ is tangent to $M$ is a local para-CR automorphism. 
\end{cor}

\section{Para-CR manifolds of finite type}
We define para-holomorphic curves as para-holomorphic mappings
\begin{align*}
F\colon \mathbb R^2 &\to \mathbb R^{2n}\\
(s,t) &\mapsto (x_1(t),\dots,x_n(t),a_1(s),\dots,a_n(s)).
\end{align*}
In other words, a para-holomorphic curve is the direct product of two curves in $\mathbb R^n$.
We call $F$ regular at $(s_0,t_0)$ if at least one of the $\frac{dx_i}{dt}(t_0)$ and at least one of the  $\frac{da_j}{ds}(s_0)$ are different from zero.

In analogy to CR geometry we define the notion of finite type $k$ at a point $p$ of an embedded para-CR manifold $M$ as the maximal order of contact of a para-holomorphic curve that is regular at $(s_0,t_0)$ and maps $(s_0,t_0)$ to $p$. 

The following proposition is analogous to the CR case:
\begin{prop}
If $M$ is a real-analytic embedded para-CR manifold then it is of finite type $k$ at $p$ if and only if there exist coordinates $(x,y,a,b)$ centred at $p$ in the ambient space such that $M$ is locally given by
\begin{equation}\label{fteq}
y=a+ \sum_{i=1}^{k-1} \gamma_i b^i x^{k-i} + o(|a|+(|x|+|b|)^k),
\end{equation}
where at least one of the $\gamma_i$ is different from zero.
\end{prop}
 
{\bf Proof.} Using the implicit function theorem we can represent $M$ in some neighbourhood of $p$ by
$$y=a + \phi(a,b,x).$$

By a coordinate change $y\mapsto y+ \phi(0,0,x)$, $a\mapsto a+ \phi(0,b,0)$ we eliminate the pure terms in $b$ and $x$ in $\phi$.

First we show that $M$ is not of finite type if all monomials in $\phi$ are divisible by $a$. Indeed,  in this case $a=0,b=s,x=t,y=0$ belongs to $M$ and therefore has contact of infinite order. 

Otherwise, $\phi$ contains a lowest order homogeneous polynomial of some degree $k$ that is not divisible by $a$. Then
$$y=a+ \sum_{i=1}^{k-1} \gamma_i b^i x^{k-i} + o(|a|+(|x|+|b|)^k)$$
where some $\gamma_i$ is different from $0$. 
Now, the para-holomorphic curve $y=a=0$, $x=t$, $b=s$ has order of contact $k$ with $M$. We show that no lower order of contact can be achieved. Indeed, $y(t)$ and $a(s)$ have to be of order $k$ at least, therefore $x(t)$ and $b(t)$ must have non-zero derivatives. This results in the order of contact $k$, which cannot be lowered. 

\hfill $\Box$
\bigskip

Notice that the highest and lowest power of $x$ that occurs in \eqref{fteq} are also invariants of the multicontact singularity $p$.

\section{Automorphisms of singular model structures}
 In this section we study the multicontact automorphisms of the singular multicontact structure $X=\partial_x$ and 
$Y=  \partial_b- (\sum_{i=1}^{k-1} i\gamma_i b^{i-1}x^{k-i})\partial_a$, $k>2$.
This is equivalent to studying the
para-CR automorphisms of the model hypersurface $S\colon y-a-P(x,b)=0$ where $P(x,b)=\sum_{i=1}^{k-1} \gamma_i b^ix^{k-i}$ for $k>2$. 

A  vector field
$$V= \alpha(a,b) \partial_a + \beta(a,b) \partial_b +\xi(x,y) \partial_x + \eta(x,y) \partial_y$$
is tangent to $S$ if $V(S)|_S=0$, namely
\begin{equation}\label{infid}
\alpha(a,b) + \beta(a,b) \sum_{i=1}^{k-1} i\gamma_i b^{i-1}x^{k-i}+ \xi(x,a+P)\sum_{i=1}^{k-1} (k-i)\gamma_i b^ix^{k-i-1} - \eta(x,a+P)=0.
\end{equation}

We assign the weight $k$ to $a$ and $y$, the weight $1$ to $b$ and $x$ and the weights $-k$ and $-1$ to
$\partial_a, \partial_y$ and  $\partial_b, \partial_x$ respectively. With respect to this grading we may consider homogeneous polynomial vector fields.
Since $S$ is homogeneous of degree $k$, it follows that the homogeneous components of an infinitesimal automorphism $V$ are infinitesimal automorphisms themselves.
  We investigate separately the homogeneous  components of an infinitesimal automorphism $V$ by considering the homogeneous components of \eqref{infid}.
\\

$\bullet$ Weight $-k$. 
In this case we have $V_{-k}= \partial_a+\partial_y$ as sole generator.
\\

$\bullet$ Weight $-k+1\le l\le 2$. Here the only terms that occur  are $\alpha_{0,\ell}b^\ell$ and
$\eta_{\ell, 0}x^\ell$. It follows immediately that $\alpha_{0\ell}=\eta_{\ell, 0}=0$.
\\

$\bullet$
Weight $-1$. In this case equation \eqref{infid} reads
$$\alpha_{0,k-1}b^{k-1}+\beta_{0,0} \sum_{i=1}^{k-1} i\gamma_i b^{i-1}x^{k-i}+ \xi_{0,0}\sum_{i=1}^{k-1} (k-i)\gamma_i b^ix^{k-i-1} - \eta_{k-1,0}x^{k-1}=0,$$
which yields
\begin{align*}
\alpha_{0,k-1}&=-\xi_{0,0}\gamma_{k-1}\\
i\beta_{0,0}\gamma_i&= -(k-i+1)\xi_{0,0}\gamma_{i-1} \quad\text{for } 2\le i\le k-1\\
\eta_{k-1,0}&=\beta_{0,0}\gamma_1
\end{align*}
From this we conclude that the only non-trivial case occurs if $S$ has the form
$$y=a+ \delta [(x+\nu b)^k-x^k-\nu^k b^k]$$
where $\delta$ and $\nu$ are determined from the relations
$$\gamma_i=\binom{k}{i}\delta \nu^i \quad\text{ for } 1\le i\le k-1.$$
The corresponding infinitesimal automorphism is a multiple of
$$V_{-1}= \partial_b- \nu \partial_x-  k\delta\nu^{k}b^{k-1} \partial_a + k\delta\nu x^{k-1}\partial_y.$$ 
After the coordinate change 
\begin{equation}\label{coordinate-change}
\begin{array}{llll}
a^*&=\frac{1}{\delta} a- \nu^k b^k & \qquad b^*&=\nu b \\
y^*&=\frac{1}{\delta} y + x^k &\qquad x^*&= x
\end{array}
\end{equation}
the equation of the surface becomes
$$y^*=a^*+ (x^*+b^*)^k$$
and the infinitesimal automorphism takes the simple form
$$V_{-1}^*=\partial_{b^*}-\partial_{x^*}.$$ 
\\

$\bullet$
Weight $0$. 
We have
$$\alpha_{1,0}a +\alpha_{0,k}b^k +\beta_{0,1}\sum_{i=1}^{k-1}i\gamma_ib^ix^{k-i} +\xi_{1,0}\sum_{i=1}^{k-1}(k-i)\gamma_ib^ix^{k-i} -\eta_{k,0}x^k - \eta_{0,1}(a+\sum_{i=1}^{k-1}\gamma_ib^ix^{k-i})=0,$$
whence
\begin{align*}
\alpha_{1,0}&= \eta_{0,1}\\
\alpha_{0,k}&=\eta_{k,0}=0\\
i\beta_{0,1}+(k-i)\xi_{1,0}&= \eta_{0,1} \qquad \forall i \text{ such that } \gamma_i\neq 0.
\end{align*}
We distinguish two cases: 
\begin{itemize}
\item[(1)] $P$ is not a monomial, that is $\gamma_i\neq 0$ and $\gamma_j\neq 0$ for two different $i$ and $j$.   Then $\beta_{0,1}=\xi_{1,0}$ and $\eta_{0,1}=k\beta_{0,1}$. The corresponding infinitesimal automorphism is
$$V_0= ka\partial_a + b\partial_b  + x\partial_x + ky \partial_y.$$  
\item[(2)]
$P= b^\iota x^{k-\iota}$ for some particular $\iota$. Then, in addition to the infinitesimal automorphism of weighted dilation from above, we have
$$V_0^\prime= (\iota-k) b\partial_b + \iota x \partial_x.$$
\end{itemize}

$\bullet$
Weight $\ell$, with $1\le\ell\le k-2$.
Using \eqref{infid} we have
\begin{multline*}
\alpha_{1,\ell}ab^\ell+\alpha_{0,k+\ell}b^{k+\ell}+ \beta_{0,\ell+1}\sum_{i=1}^{k-1} i\gamma_ib^{i+\ell}x^{k-i} + \xi_{\ell+1,0} \sum_{j=1}^{k-1} (k-j)\gamma_jb^{j}x^{k+\ell-j} \\- \eta_{k+\ell,0}x^{k+\ell}-\eta_{\ell,1}x^\ell(a+\sum_{i=1}^{k-1}\gamma_ib^ix^{k-i}) =0.
\end{multline*}
It readily follows that 
$$\alpha_{1,\ell}=\eta_{\ell,1}=\alpha_{0,k+\ell}=\eta_{k+\ell,0}=0.$$
We show that $\beta_{0,\ell+1}=\xi_{\ell+1,0}=0$. Suppose they do not vanish. Since the powers $b^j$ with $1\le j\le \ell$ appear only in the sum $\sum_{j=1}^{k-1} (k-j)\gamma_jb^{j}x^{k+\ell-j}$, we immediately have $\gamma_j=0$ for all such $j$. Now iteration of the same argument shows that all $\gamma_j$ vanish, which is not possible. The contradiction shows that $\beta_{0,\ell+1}=\xi_{\ell+1,0}=0$.
\\

$\bullet$
Weight $k-1$.
Equation  \eqref{infid} yields
\begin{multline}\label{eq6}
\alpha_{1,k-1}ab^{k-1}+\alpha_{0,2k-1}b^{2k-1}+ \beta_{0,k}\sum_{i=1}^{k-1} i\gamma_ib^{i+k-1}x^{k-i} + \beta_{1,0}a\sum_{i=1}^{k-1} i\gamma_ib^{i-1}x^{k-i}+\\ \xi_{k,0} \sum_{j=1}^{k-1} (k-j)\gamma_jb^{j}x^{2k-j-1} +  \xi_{0,1}(a+\sum_{i=1}^{k-1}\gamma_ib^ix^{k-i}) \sum_{j=1}^{k-1} (k-j)\gamma_jb^{j}x^{k-j-1}\\- \eta_{2k-1,0}x^{2k-1}-\eta_{k-1,1}x^{k-1}(a+\sum_{i=1}^{k-1}\gamma_ib^ix^{k-i}) =0
\end{multline}
It follows immediately that $\alpha_{1,k-1}= \xi_{0,1}\gamma_{k-1}$ and
$\alpha_{0,2k-1}=\eta_{2k-1,0}=0$.
Inspecting the coefficient of $ax^{k-1}$ we get $\beta_{1,0} \gamma_1=\eta_{k-1,1}$.
Similarly, looking at $ab^jx^{k-j-1}$ we get $(j+1)\beta_{1,0}\gamma_{j+1}=(k-j)\xi_{0,1}\gamma_j$ for $1\le j\le k-1$.
Therefore, we have that either $\beta_{1,0}=\xi_{0,1}=0$ (which entails that all other coefficients vanish) or $\gamma_j=\binom{k}{j}\delta \nu^j$ with
$$\nu=\frac{\xi_{0,1}}{\beta_{1,0}}, \qquad \delta=\frac{\eta_{k-1,1}}{k\xi_{0,1}}.$$
%
%
After the coordinate change in \eqref{coordinate-change} we may assume that $S$ has the form
$$y=a+(x+b)^k.$$
In this case $\partial_b-\partial_x$ is an infinitesimal automorphism of weight $-1$. Its commutator with an infinitesimal automorphism $V$ of weight $k-1$ would have weight $k-2$, which does not exist if $k>2$. Therefore $[\partial_b-\partial_x,V]=0$, which implies that $V$ must be a multiple of $y\partial_x$. However, this is not an infinitesimal automorphism of $S$.
\\

$\bullet$
Weight $k$. An infinitesimal automorphism of weight $k$ has the form
\begin{multline*}
V=(\alpha_{2,0}a^2+\alpha_{1,k}ab^k+ \alpha_{0,2k}b^{2k})\partial_a+ (\beta_{1,1}ab+\beta_{0,k+1}b^{k+1})\partial_b +\\ + (\xi_{1,1}xy+ \xi_{k+1,0}x^{k+1})\partial_x + (\eta_{0,2}y^2+\eta_{k,1}x^ky+\eta_{2k,0}x^{2k})\partial_y.
\end{multline*}
Checking the power $2k$ and the next highest power in $b$ and $x$ in 
$V(S)|_S=0$,
we conclude that 
$$\alpha_{0,2k}=\beta_{0,k+1}=\xi_{k+1,0}=\eta_{2k,0}=0.$$
Now, $[\partial_y+\partial_a,V]$ must be an infinitesimal automorphism of weight $0$ and therefore is a linear combination of $a\partial_a$, $b\partial_b$, $x\partial_x$, $y\partial_y$.  This implies that
$$\alpha_{1k}=\eta_{k,1}=0.$$
Suppose that $P$ is not a monomial. Then  
$$[\partial_y+\partial_a,V]= s (x\partial_x+b\partial_b+ k y\partial_y+ k a\partial_a)$$
for some $s\in \mathbb R$
and therefore
$$V=s(\frac{k}{2}a^2\partial_a+ ab \partial_b+ xy\partial_x+ \frac{k}{2}y^2\partial_y).$$
The restriction of $V(S)|_{S=0}$ to $a=0$ yields
$$(\sum_{i=1}^{k-1} \gamma_i b^i x^{k-i})(\sum_{i=1}^{k-1}(\frac{k}{2}-i) \gamma_i b^i x^{k-i})=0.$$
This means that $i=\frac{k}{2}$ which contradicts our assumption. Hence, without loss of generality, $P=b^\iota x^{k-\iota}$ and 
$$[\partial_y+\partial_a,V]= ska\partial_a+(s+t(\iota-k))b\partial_b+(s+\iota t)x\partial_x+sky\partial_y.$$
It follows
$$V=\frac{sk}{2}a^2\partial_a+(s+t(\iota-k))ab\partial_b+(s+\iota t)xy\partial_x+\frac{sk}{2}y^2\partial_y.$$
The restriction of $V(S)|_{S=0}$ to $a=0$ yields
$$(k-\iota)(s+\iota t)-\frac{sk}{2}=0,$$
whence 
$$t=\frac{\iota-\frac{k}{2}}{\iota(k-\iota)}s.$$
We conclude that a generator of an infinitesimal automorphism of weight $k$ is
$$V_k=a^2\partial_a+\frac{1}{\iota}ab\partial_b+\frac{1}{k-\iota}xy\partial_x+y^2\partial_y.$$
\\

$\bullet$
Weights $k+\ell$, with $l>0$.
 Let $V$ be a vector field of homogenous weight $k+\ell$, namely
$$V=\sum_{ki+j=2k+\ell} \alpha_{i,j}a^ib^j \partial_a+ \sum_{ki+j=k+\ell+1} \beta_{i,j}a^ib^j \partial_b+ \sum_{i+kj=k+\ell+1} \xi_{i,j}x^iy^j \partial_x+ \sum_{i+kj=2k+\ell} \eta_{i,j}x^iy^j \partial_y.$$
The commutator $[\partial_a+\partial_y,V]=0$, except for the case when $\ell$ is a multiple of $k$ and $P$ is a monomial, which we will treat separately.
It readily follows that $V$ does not depend on $a$ and $y$, i.e.
 $$V= \alpha_{0,2k+\ell}b^{2k+\ell} \partial_a+ \beta_{0,k+\ell+1}b^{k+\ell+1} \partial_b+ \xi_{k+\ell+1,0}x^{k+\ell+1} \partial_x+ \eta_{2k+\ell,0}x^{2k+\ell}\partial_y.$$
Inspecting the highest order terms in $V(S)|_S=0$ yields that $V=0$.
Now, if $\ell=2k$ and $P=b^mx^n$, then
$$[\partial_a+\partial_y,V]= t (a^2\partial_a+\frac{1}{m}ab\partial_b+\frac{1}{n}xy\partial_x+y^2\partial_y)$$
for some $t\in \mathbb R$
and therefore
$$V= \frac{t}{3}a^3\partial_a+\frac{t}{2}a^2b\partial_b+\frac{t}{2}xy^2\partial_x+\frac{t}{3}y^3\partial_y.$$
The restriction of $V(S)|_S=0$ to $a=0$ yields
$$\frac{t}{6}b^{3m}x^{3n}=0,$$
whence $t=0$. The previous argument shows that, also in the case when $P$ is a monomial, there are no infinitesimal automorphisms of weight higher than $k$.
\\
We summarize the conclusions from the above computations in the following statement.
\begin{theorem}\label{main}
Let $S\colon y-a-P(x,b)=0$ be a hypersurface with $P(x,b)=\sum_{i=1}^{k-1} \gamma_ib^ix^{k-i}$, where $k>2$ and $\gamma_i\in \mathbb R$, with $\gamma_i\neq 0$ for some $i\le k-1$. For the space $\mathcal A$ of multicontact infinitesimal automorphisms we have three cases.
\begin{itemize}
\item[(i)] If $P(x,b)=b^\iota x^{k-\iota}$, then $\mathcal A ={\mathbb R} V_{-k}+ {\mathbb R} V_{0}+ {\mathbb R} V_{0}^\prime + {\mathbb R} V_{k}$. As a Lie algebra, $\mathcal A$ is isomorphic to
$ {\mathfrak sl}(2,\mathbb R)\oplus \mathbb R$.
\item[(ii)] If $P(x,b)=a+(x+b)^k$, then $\mathcal A = {\mathbb R} V_{-k}+{\mathbb R} V_{-1}+ {\mathbb R} V_{0}$. The space $\mathcal A$ is a solvable Lie algebra isomorphic to the $\mathfrak n \oplus \mathfrak a$ part of the Iwasawa decoposition of ${\mathfrak so}(1,5)$.
\item[(iii)] If $P(x,b)$ is not as in (i) or (ii), then $\mathcal A= {\mathbb R} V_{-k}+{\mathbb R} V_{0}$. As a Lie algebra, $\mathcal A$ is isomorphic to
 the  $\mathfrak n \oplus \mathfrak a$ part of the Iwasawa decoposition of ${\mathfrak sl}(2,\mathbb R)$.
\end{itemize}
\end{theorem}
We may integrate the vector fields that give the infinitesimal automorphisms to obtain $1$-parameter groups of automorphisms. For (iii) in the theorem above, we have
\begin{align}
{\rm Exp}(t{V_{-k}}): (x,y,a,b)&\to (x,y+t,a+t,b),\quad t\in \mathbb R \label{14}\\
{\rm Exp}(\lambda{V_{0}}): (x,y,a,b)&\to (\lambda x,\lambda^k y,\lambda^k a,\lambda b),\quad \lambda >0 \label{15}
\end{align}
In the case (i) we add to \eqref{14} and \eqref{15} the flows
\begin{align*}
{\rm Exp}(\lambda V_0^\prime):(x,y,a,b)&\to \left(\lambda^\iota x,y,a,\frac{b}{\lambda^{k-\iota}}\right),\quad \lambda >0\\
{\rm Exp}(t{V_{k}}): (x,y,a,b)&\to \left(\frac{x}{\sqrt[k-\iota]{1-ty}},\frac{y}{(1-ty)},\frac{a}{1-ta},\frac{b}{\sqrt[\iota]{1-ta}}\right).
\end{align*}
Finally, in the case (ii) we have  \eqref{14} and \eqref{15} and the following flow that we write in the coordinates \eqref{coordinate-change}
$$
{\rm Exp}(tV_{-1}):(x,y,a,b)\to (x-t,y-x^k+2(x-t)^k,a+b^k-2(b+t)^k,b+t).
$$

Below we compute the multicontact automorphisms that are not exponential images of infinitesimal automorphisms. First we show that they have to be linear. 
Let $S\colon y=a+ P(x,b)$ with $P=a+ \sum_{j=1}^{k-1} \gamma_j b^j x^{k-j}$ and the distinguished direction fields $X=\partial_x+ P_x \partial_y$ and $Y=\partial_b- P_b \partial_a$, as above. The induced multicontact structure has a singularity at the zero locus of 
$[X,Y]=-P_{xb} (\partial_a+\partial_y)$, i.e. for $P_{xb}=0$. The polynomial $P_{xb}$ is homogeneous of degree $k-2>0$ and factorises into a product of irreducible quadratic and linear factors. Therefore, the zero locus is either a point (if all factors are irreducible quadratic polynomials), a line (if $P_{xb}=\#(\mu x+ \nu b)^{k-2}$) or a pencil of lines passing through $0$ (if $P_{xb}$ contains at least two different linear factors.) In the first and the last case the origin in the $xb$-plane is preserved by any multicontact automorphism. In the second case one can choose coordinates such that
$P=(x+b)^k-x^k-b^k$. In this case the line $x+b=0$ is preserved by any multicontact automorphism. 

\begin{lem} Any automorphism $\Phi$ is a composition $\Phi=\Phi_0\circ \Phi_1$, where $\Phi_0(0)=0$ and $\Phi_1$ is of the form
$$y\mapsto y+ t, \quad a\mapsto a+t$$
if $P_{xb}$ is not a power of a linear term, and
$$y\mapsto y+ t, \quad a\mapsto a+t, \quad x\mapsto x+s, \quad b \mapsto b-s$$
if $P_{xb}=\# (x+b)^{k-2}$.
\end{lem}

{\bf Proof.} Let $x\mapsto \Xi(x,y)$, $y\mapsto H(x,y)$, $a\mapsto A(a,b)$, $b\mapsto B(a,b)$ be an automorphism of $S$. If $P_{xb}$ is not a power of a linear term, we have $B(0,0)=0$ and $H(0,0)=0$.
Inspecting the zero order term in $H(x,a+P)=A(a,b) + P(B(a,b),\Xi(x,a+P))$ we find $H(0,0)=A(0,0)=t$.
If we split $\Phi$ into $\Phi_0\circ \Phi_1$ then $\Phi_0$ clearly preserves the origin.

Consider the case  $P_{xb}=\# (x+b)^{k-2}$. Inspecting the zero order term in $H(x,a+P)=A(a,b) + P(B(a,b),\Xi(x,a+P))$ we find $H(0,0)=A(0,0)+ P(B(0,0),\Xi(0,0))$. Since $x+b=0$ is preserved, it follows $\Xi(0,0)=-B(0,0)=s$ and $P(B(0,0),\Xi(0,0))=0$, hence again $H(0,0)=A(0,0)=t$. For $\Phi=\Phi_0\circ \Phi_1$, we conclude that $\Phi_0$ preserves the origin. \hfill $\Box$

\begin{lem} Any automorphism $\Phi$ with $\Phi(0)=0$ splits into $\Phi=\Phi_3\circ\Phi_2$ where $\Phi_2$ is of the form $x\mapsto t x$, $b\mapsto s b$, $y\mapsto \tau y$, $a\mapsto \tau a$, where 
$\tau= s^i t^j$ for all $i,j=k-i$ such that $\gamma_i\neq 0$.
\end{lem}

{\bf Proof.} Let $\Phi_2$ be the linear part of $\Phi$, i.e., $\Phi= \Phi_2+ O(|x|+|b|+|y|+|a|)$. Then 
$\Phi_2$ is an automorphism itself. Inspecting the coefficients of the linear terms yields the desired result. \hfill $\Box$

\begin{lem} Any real-analytic automorphism $\Phi_3$ is of the form $\Phi_3=\exp V$ for some infinitesimal automorphism $V$.
\end{lem}

{\bf Proof.} This  is similar to Lemma 3.7 in \cite{Sch98}. \hfill $\Box$

\begin{prop}
The group of discrete multicontact automorphisms of $S$ is either $\mathbb Z^2$ and generated by 
$$x\mapsto -x, \qquad b\mapsto -b, \qquad a\mapsto (-1)^k a, \qquad  y\mapsto (-1)^k y$$
or $\mathbb Z^2\times \mathbb Z^2$ with second generator
$$x\mapsto x, \qquad b\mapsto -b, \qquad a\mapsto (-1)^i a, \qquad  y\mapsto (-1)^i y,$$ if
all $i$ with $\gamma_i\neq 0$ are of the same parity.
\end{prop}
{\bf Proof.} After combining a linear automorphism  with suitable one-parametric families we may restrict to $x\mapsto \pm x$ and  $b\mapsto \pm b$. For $x\mapsto -x$, $y\mapsto -y$ we always have an automorphism with $a\mapsto (-1)^k a$  and $ y\mapsto (-1)^k y$. For $x\mapsto x$, $b\mapsto -b$ we only get an automorphism if $(-1)^i$ is the same for all $i$ with $\gamma_i\neq 0$. The remaining case  $x\mapsto -x$, $b\mapsto b$  is a composite of the previous ones.    \hfill $\Box$

\begin{remark}
We notice that the hypersurface 
$$y=a+P(b,x)=a+ \sum_{j=1}^{k} \gamma_j b^j x^{k-j},$$
can be seen as the manifold of solutions of the ODE
$$
y^{(k)}(x)=0
$$
with initial conditions $y(0)=a$ and  $y^{(j)}(0)= j!\gamma_{k-j}b^{k-j}$ for $j=1,\dots,k-1$.


Case (i) of Theorem \ref{main} is the hypersurface
$$y=a+bx^\ell,$$
which can be considered as the manifold of solutions of the ODE
$$y^{(\ell+1)}=0$$
with initial conditions $y(0)=a$, $y'(0)=\cdots=y^{(\ell-1)}(0)=0$, $y^\ell (0)=b$.
Alternatively, it is the manifold of solutions of the singular 2nd order ODE $xy''-(\ell-1) y'=0$.
\\

Case (ii)  of Theorem \ref{main} is the hypersurface $y=a+(x+b)^k$, which can be viewed  as the solution of
$$y''= s (y')^t$$
with $s=(k-1)k^{\frac{1}{k-1}}$, $t=\frac{k-2}{k-1}$.

\end{remark}

%
%

\end{document}